# Some Observations on the 3x+1 Problem

Dhananjay P. Mehendale
Sir Parashurambhau College, Tilak Road, Pune-411030, India

Abstract

We present some interesting observations on the 3x+1 problem. We propose a new algorithm which eliminates certain steps while we check the action of 3x+1 procedure on a number. Also, we propose a reason why many numbers follow a similar pattern during execution of 3x+1 algorithm. We advocate towards the end (a heuristic argument) that the 3x+1 conjecture is more likely to be true than false.

**1. Introduction:** Perhaps the most fascinating problem in number theory is the so called 3x+1 problem. This problem is also known as Collatz's problem, Kakutani's problem, Syracuse problem, Ulam's problem, and Hasse's algorithm [1]. The problem **can be** stated as follows:

Let $T: Z \to Z$ be defined by

$$T(x) = \frac{x}{2^k} \quad \text{if } x \equiv 0 \bmod(2), \text{ and } k \text{ is the}$$

largest power of 2 that divides $x$.

$$T(x) = \frac{(3x+1)}{2^k} \quad \text{if } x \equiv 1 \bmod(2), \text{ and } k \text{ is the}$$

largest power of 2 that divides $(3x+1)$.

It is **conjectured** that if $x \in N$ then the trajectory $x, T(x), T^2(x), \cdots$ eventually reaches (converges to) 1 in finitely many steps.

We now proceed to discuss

**2. Some Interesting Observations:**

**Observation (1):** (a) For an odd number $x$ of type $(4k+1)$ we have $T(x) < x$.

(b) For an odd number $x$ of type $(4k+3)$ we have $T(x) > x$.

**Proof of (a):** $T(x) = \dfrac{3x+1}{2^u} = \dfrac{12k+4}{2^u} = \dfrac{4(3k+1)}{2^u}$, hence divisible at least by 4 (i.e. $u$ can be at least equal to 2).

**Proof of (b):** $T(x) = \dfrac{3x+1}{2^u} = \dfrac{12k+10}{2^u} = \dfrac{2(6k+5)}{2^u}$, hence divisible at most by 2 (i.e. $u$ can be at most equal to 1).

**Observation (2):** Every odd integer $y$ can be uniquely expressed as
$$y = 2^n x + 2^{(n-1)} - 1,$$
i.e. for every odd integer $y$ there exist unique $x$, $n$ such that the above equation holds.

**Proof:** Find the unique largest power of 2 that divides $(y+1)$, say $(n-1)$, and let $\dfrac{(y+1)}{2^{(n-1)}} = 2x+1$. Hence etc.

**Observation 3:** There are numbers which rise $n$ number of times where $n$ is any arbitrary positive integer without a single fall under the action of $3x + 1$ algorithm.

**Proof:** Consider the following number
$$y = 2^{(n+2)} x + 2^{(n+1)} - 1, \text{ then}$$
$$T(y) = 2^{(n+1)}(3x+1) + 2^n - 1,$$

$$T^2(y) = 2^n(3^2 x + 3^1 + 3^0) + 2^{(n-1)} - 1,$$
$$\vdots$$
$$T^n(y) = 2^2(3^n x + \dfrac{(3^n - 1)}{2}) + 2^1 - 1 = y_1.$$

After this stage we get
$$T^{(n+1)}(y) = T(y_1).$$

Note that $T^j(y), j = 1, 2, (n-1)$ is an odd number of type $(4k+3)$ while $T^n(y)$ is an odd number of type $(4k+1)$. Therefore, it is clear that
$$y < T^1(y) < T^2(y) < \cdots < T^n(y), \text{ while}$$
$$T^{(n+1)}(y) < T^n(y).$$

This observation offers us a faster algorithm by which we can omit certain steps while we process a positive odd integer by $3x + 1$ algorithm as follows:

**A Faster $3x + 1$ Algorithm:**
   (1) Express the given positive odd integer $y$ uniquely as

$$y = 2^n x + 2^{(n-1)} - 1$$

(2) Find directly

$$T^{(n-2)}(y) = 2^2(3^{(n-2)}x + \frac{3^{(n-2)} - 1}{2}) + 2^1 - 1 = y_1, \text{ say}$$

(3) Find $T(y_1) = y_2$, say and if $y_2 \neq 1$ then set $y \leftarrow y_2$ and go to step (1).

**Example 1:** Let $y = 1023$, then $T^8(y) = y_1 = 39365$ and
$T(y_1) = y_2 = 7381$ (= a new $y$ to start with). Thus, 8 executions are **avoided!!**

**Observation 4:** If $y \to 1$ then $2^k y \to 1$ for all $k \geq 0$.

**Proof:** Obvious.

**Observation 5:** If $y \to 1$ then $z = 4^n y + (\frac{4^n - 1}{3}) \to 1$
for all $n = 1, 2, \ldots$ More precisely, if $T(y) = y_0$ then $T(z) = y_0$.

**Proof:** We proceed by induction on $n$.

**Step (1):** Let $n = 1$. Therefore, we have $z = 4y + 1$.
Now, $T(y) = y_0 = \frac{(3y+1)}{2^j}$ where $j$ is the largest power of 2 that
divides $(3y + 1)$. Also, $3z + 1 = 3(4y + 1) + 1 = 4(3y + 1)$, therefore the
largest power that divides $(3z + 1)$ will be $2^{(j+2)}$ and
$T(z) = y_0 = \frac{(3z+1)}{2^{(j+2)}}$.

**Step (2):** We assume by induction the result holds for $n = k$ and
proceed to show it for $n = (k + 1)$. Let $z = 4^k y + (\frac{4^k - 1}{3})$ and
$T(z) = y_0$. We have to see that when $z' = 4^{(k+1)} y + (\frac{4^{(k+1)} - 1}{3})$ we
still have $T(z') = y_0$. Note that $z' = 4z + 1$ and when $T(z) = y_0$ then
also $T(z') = y_0$ by step (1).

We now proceed to see why many numbers follow a similar pattern during execution of $3x+1$ algorithm.

**Numbers with Similar Pattern:** From observations 4 and 5 suppose a number $x_0$ has the following sequence under the action of $3x + 1$ algorithm:

$$x_0 \to x_1 \to x_2 \to \cdots \to x_j \to \cdots \to x_n \to 1, \text{ then}$$

(1) The numbers $2^k x_j$ will have same pattern for all $k, j$.

(2) The numbers $4^k x_j + (\dfrac{4^k - 1}{3})$ will have same pattern for all $k, j$.

**Example 2:** Count the numbers $\leq 1000$ which have maximum value 9232.

**Solution:** Using (1) and (2) we count the numbers as follows:

(1) Find the smallest number $x_0^1$ say such that $x_0^1 \leq 1000$.

Let $x_0^1 \to x_1^1 \to \cdots \to x_j^1 \to \cdots \to x_k^1 \to 9232$.

(2) Consider the set of numbers $\{2^m (4^n x_j^1 + \left(\dfrac{4^n - 1}{3}\right))\}$,

for all nonnegative integers $m, n$ and $0 \leq j \leq k$, such that all these numbers are $\leq 1000$. It is clear that all these numbers when processed under $3x+1$ algorithm will have 9232 as maximum.

(3) Find all the numbers $z \leq 1000$ which **fall** to $x_j^1, 0 \leq j \leq k$, under the $3x + 1$ algorithm. Consider all distinct numbers $z_i$ on the paths from $z$ to $x_j^1$ and find as above their multiples $\{2^m (4^n z_i + \left(\dfrac{4^n - 1}{3}\right))\}$, for all nonnegative integers $m, n$ and $0 \leq i \leq k$, such that all these numbers are $\leq 1000$. Again, all these numbers when processed under $3x+1$ algorithm will have 9232 as maximum.

(4) Consider next smallest number other than the numbers considered above having maximum 9232, say $x_0^2$, etc. etc. and go to (1).

(5) Continue till all the numbers $\leq 1000$ that go to 9232 are considered in some set.

One can easily check that there are **in all 350 numbers** $\leq 1000$ having **maximum 9232** under the $3x+1$ algorithm.

**Observation 6:** Solving the 3x+1 conjecture is equivalent to showing that every positive odd integer $x$ has a **representation** as:

$$x = \frac{2^{n_{(k+1)}} - \{3^k 2^0 + 3^{(k-1)} 2^{n_1} + \cdots + 3^0 2^{n_k}\}}{3^{(k+1)}} \quad \rightarrow (1)$$

where the integral indices satisfy $0 < n_1 < n_2 < \cdots < n_{(k+1)}$, $n_{(k+1)} \neq n_k + 2$.

**Proof:** Simple.

**Example 3:** $7 = \dfrac{2^{11} - \{3^4 2^0 + 3^3 2^1 + 3^2 2^2 + 3^1 2^4 + 3^0 2^7\}}{3^5}$

**Observation 7:** If we consider **all** the solutions (integral as well as rational) of equation (1) we get **all possible** $x$ that go to 1 in (k+1) steps. Consider following sets of **integral** and **rational** solutions:

Let $U_1 = \left\{ x / x = \dfrac{2^k - 3^0 2^0}{3^1}, k = 1, 2, \cdots \right\}$

$U_2 = \left\{ x / x = \dfrac{2^{k_2} - \{3^1 2^0 + 3^0 2^{k_1}\}}{3^2}, 0 < k_1 < k_2, k_1, k_2 = 1, 2, \cdots \right\}$, etc.

**Observation 8:** Solving 3x+1 conjecture is equivalent to showing that

$$Z_{Odd} \subset \bigcup_{j=1}^{\infty} U_j$$

**Why 3x+1 conjecture is more likely to be true than false?** Suppose a number $x_0$ requires $k$ steps to become 1. One can see (by numerically tackling some examples) that if one finds **closest possible** $x_1^i, x_2^i$ belonging to $U_i$, $1 \leq i \leq (k-1)$ such that $x_1^i < x_0 < x_2^i$ then one observes that one moves **closer and closer** to $x_0$ as $i$ is increased from 1 to (k-1) (as one gets more parameters for maneuvering) and finally at k-th step one gets a solution in $U_k$ which is equal to $x_0$.